\def\squareforqed{\hbox{\rlap{$\sqcap$}$\sqcup$}}
\def\qed{\ifmmode\squareforqed\else{\unskip\nobreak\hfil
\penalty50\hskip1em\null\nobreak\hfil\squareforqed
\parfillskip=0pt\finalhyphendemerits=0\endgraf}\fi}
\def\bbbr{{\rm I\!R}} 
\documentclass[12pt]{article}
\def\d{\displaystyle}
\begin{document}
\newtheorem{theorem}{Theorem}
\newtheorem{col}{Corollary}
\newtheorem{lem}{Lemma}
\newtheorem{prop}{Proposition}
\newtheorem{rem}{Remark}
\title{Positivity of Tur\'an determinants for orthogonal polynomials}
\author{Ryszard Szwarc\thanks{\noindent This work has been partially supported
by KBN (Poland) under grant 2 P03A 030 09.}}
\date{}
\maketitle
\begin{abstract}
The orthogonal polynomials $p_n$ satisfy Tur\'an's inequality if
$p_n^2(x)-p_{n-1}(x)p_{n+1}(x)\ge 0$ for $n\ge 1$ and for all $x$ in the
interval of orthogonality.  We give general criteria for orthogonal
polynomials  to satisfy Tur\'an's inequality.  This yields the
known results for classical orthogonal polynomials as well as new results,
for example, for the $q$--ultraspherical polynomials.
\end{abstract}

\footnotetext{\noindent 1991 {\it Mathematics Subject Classification}.
Primary  42C05, 47B39}
\footnotetext{\noindent
{\it Key words and phrases}: orthogonal polynomials, Tur\'an's
inequality, recurrence
formula.}
\section{Introduction}
In the 1940's, while studying the zeros of Legendre polynomials $P_n(x),$
Tur\'an \cite{tu} discovered that
\begin{equation}\label{turan}
P_n^2(x)-P_{n-1}(x)P_{n+1}(x) \ge 0, \qquad -1\le x\le 1
\end{equation}
with equality only for $x=\pm 1.$
Szeg\"o \cite{szego} gave four different proofs of (\ref{turan}).
Shortly after that, analogous results were obtained for other classical orthogonal
polynomials such as  ultraspherical polynomials \cite{sk,szasz},
Laguerre and  Hermite polynomials
\cite{mn}, and  Bessel functions \cite{sk,szasz}.

In \cite{ks} Karlin and Szeg\"o raised the question of determining the range of parameters
$(\alpha,\beta)$ for which  (\ref{turan}) holds for Jacobi polynomials
of order $(\alpha,\beta);$ i.e. denoting $R_n^{(\alpha ,\beta )}(x)=
P_n^{(\alpha ,\beta )}(x)/P_n^{(\alpha ,\beta )}(1),$
\begin{equation}\label{jacobi}
[R_n^{(\alpha ,\beta )}(x)]^2-R_{n-1}^{(\alpha ,\beta
)}(x)R_{n+1}^{(\alpha ,\beta )}(x)\ge 0, \qquad -1\le x \le 1.
\end{equation}
In 1962 Szeg\"o \cite{szego2} proved (\ref{jacobi}) for $\beta\ge
|\alpha|,\ \alpha>-1.$ In a series of two papers \cite{g1,g2} Gasper
extended Szeg\"o's result by showing that (\ref{jacobi}) holds if and only
if $\beta \ge \alpha >-1.$

More recently, attention has also turned to the q-analogues of the
classical polynomials \cite{bi}.

All the results mentioned above were proved using differential
equations, that the
classical orthogonal polynomials satisfy.  Therefore the methods cannot
be used to extend (\ref{turan}) to more general orthogonal polynomials.
In 1970 Askey \cite[Thm. 3]{a1} gave a general criterion for monic symmetric
orthogonal polynomials to satisfy the Tur\'an type inequality on the
entire real line.  His
result, however, does not imply (\ref{turan}) for the Legendre polynomials
because the latter are not monic in the standard normalization, and they
do not satisfy Askey's assumptions
in
the monic normalization.
In this paper we give general criteria for orthogonal polynomials
implying (\ref{turan}) holds for $x$ in the support of corresponding
orthogonality measure.  The assumptions are stated in terms of the
coefficients of the recurrence relation that the orthogonal polynomials
satisfy.  They admit a very simple form in the case of symmetric
orthogonal polynomials; i.e.  the case
$p_n(-x)=(-1)^np_n(x).$ In particular, the results  apply to all
the ultraspherical polynomials, giving yet another proof of Tur\'an's
inequality for the Legendre polynomials.

It turns out that the way we normalize the polynomials is essential for
the Tur\'an inequality to hold. The results concerning the classical
orthogonal polynomials used the normalization at  one
endpoint of the interval of orthogonality, e.g.  at $x=1$ for the
Jacobi polynomials and at $x=0$ for the Laguerre polynomials.
We will also use this normalization and will show that this choice is
optimal (Proposition \ref{prop}).
However, the recurrence relation for the polynomials normalized in this
way may not be available explicitly. This is the case of the
$q$--ultraspherical
polynomials. We give a way of overcoming this obstacle (Corollary 1).
In particular, we prove the Tur\'an inequality for all q-ultraspherical
polynomials with $q>0.$ These polynomials have been studied by Bustoz
and Ismail \cite{bi} but with a normalization other than at $x=1.$
The same method is applied to the symmetric Pollaczek polynomials,
studied in \cite{bi1}, again with different normalization.

In  Section 6  we prove results for nonsymmetric orthogonal
polynomials (Thm.  4).  The assumptions again are given in terms of the
coefficients in a three term recurrence relation but they are much more
involved.

In Section 7 we state results concerning polynomials orthogonal on the
positive half axis.  In particular they can be applied to the Laguerre
polynomials of any order $\alpha .$

\section{Basic formulas.}
Let $p_n$ be polynomials orthogonal with respect to a
probability measure on $\bbbr.$
The expressions
\begin{equation}\label{det}
\Delta_n (x)= p_n^2(x)-p_{n-1}(x)p_{n+1}(x) \quad
n=0,\,1,\,\ldots ,
\end{equation}
are called the Tur\'an determinants.
Our goal is to give conditions implying the nonnegativity of $\Delta_n
(x)$ for $x$ in the support of the orthogonality measure.

The first problem we encounter is that the orthogonality
 determines the
polynomials $p_n$ up to a nonzero  multiple.
The sign of $\Delta_n(x)$ may change if we multiply
each $p_n$ by  different nonzero constants.
We will normalize the polynomials $p_n$ to obtain the sharpest
results possible. Namely, we will assume
that
$$p_n(a)=1$$  at a point $a$ in the support of the
orthogonality measure. In this way the Tur\'an determinant
vanishes at $x=a.$

Our main interest is focused on the case when the orthogonality
measure is supported in a bounded interval.
By an affine change of variables we can assume
that this interval is $[-1,1].$ In that case we set $a=1.$
Since the polynomials $p_n$ do not change sign in the interval
$[1,+\infty)$ they have positive leading coefficients.

Assume that the polynomials $p_n$ are orthogonal, with positive leading
coefficients and $p_n(1)=1.$
Then they satisfy
 the three term recurrence relation
\begin{equation}\label{rec}
xp_n(x)=\gamma_np_{n+1}(x)+\beta_np_n(x) +
\alpha_np_{n-1}(x)\quad n=0,\,1,\,\ldots,
\end{equation}
with initial conditions
$p_{-1}=0,$ $p_0=1,$ where $\alpha_n,$
$\beta_n,$ and $\gamma_n$ are given sequences of real valued
coefficients such that
$$\alpha_0=0, \quad \alpha_{n+1}>0,\quad \gamma_n >0 \quad {\rm for} \
n=0,\,1,\,\ldots .$$
Plugging $x=1$ into  (\ref{rec}) gives
\begin{equation}\label{norm}
\alpha_n + \beta_n + \gamma_n =1\qquad  n=0,\,1,\,\ldots .
\end{equation}

\begin{prop}\label{prop}
Let the polynomials $p_n$ satisfy (\ref{rec}) and (\ref{norm}).
Then
\begin{eqnarray}
\gamma_n\Delta_n &\hspace{-8pt}=\hspace{-8pt}&  \gamma_np_n^2+\alpha_np_{n-1}^2 -
(x-\beta_n)p_{n-1}p_n,\label{1}\\
\gamma_n\Delta_n &\hspace{-8pt}=\hspace{-8pt}&
(p_{n-1}-p_{n})[(\gamma_{n-1}-\gamma_n)p_n+
(\alpha_n-\alpha_{n-1})p_{n-1}] +\alpha_{n-1}\Delta_{n-1},\label{2}\\
\gamma_n\Delta_n &\hspace{-8pt}=\hspace{-8pt}&
(p_{n}-p_{n-1})(\gamma_np_n-\alpha_np_{n-1})+ (1-x)p_{n-1}p_n,\label{3}
\end{eqnarray}
for $n=1,\,2,\,\ldots\,.$
\end{prop}

\noindent{\it Proof.}
By (\ref{rec}) we get
\begin{eqnarray*}
\gamma_n\Delta_n&=&\gamma_np_n^2-\gamma_np_{n-1}[(x-\beta_n)-\alpha_np_{n-1}]\\
&=&\gamma_np_n^2+\alpha_np_{n-1}^2 -
(x-\beta_n)p_{n-1}p_n\\
&=&\gamma_np_n^2+\alpha_np_{n-1}^2 -
(\beta_{n-1}-\beta_n)p_{n-1}p_n -(x-\beta_{n-1})p_{n-1}p_n.
\end{eqnarray*}
Now applying (\ref{rec}), with $n$ replaced by $n-1,$ to the last term
yields
$$
\gamma_n\Delta_n=(\gamma_n-\gamma_{n-1})p_n^2 +
(\alpha_n-\alpha_{n-1})p_{n-1}^2
-(\beta_{n-1}-\beta_n)p_{n-1}p_n + \alpha_{n-1}\Delta_{n-1}.$$
The use of
$$ \beta_{n-1}-\beta_n= (\gamma_n-\gamma_{n-1}) +
(\alpha_n-\alpha_{n-1})$$  concludes the proof of (\ref{2}).
In order to get (\ref{3}) replace $\beta_n$ with
$1-\alpha_n-\gamma_n$ in  (\ref{1}).
\qed

\section{Symmetric polynomials}
We will consider first the symmetric orthogonal polynomials,
i.e. the orthogonal polynomials satisfying
\begin{equation}\label{-x}
 p_n(-x)=(-1)^np_n(x).
\end{equation}
\begin{theorem}
Let the polynomials $p_n$ satisfy
\begin{equation}\label{recsym}
xp_n(x)=\gamma_np_{n+1}(x) + \alpha_np_{n-1}(x)\qquad
n=0,\,1,\,\ldots .
\end{equation}
with $p_{-1}=0,$ $p_0=1,$ where $\alpha_0=0,$ $\alpha_{n+1}>0,\,
\gamma_n>0,$ and
$$\alpha_n+\gamma_n=a \qquad n=0,\,1,\,\ldots .$$
Assume that either (i) or (ii) is satisfied where
\begin{enumerate}
\item[(i)]
$\alpha_n$ is nondecreasing and $\alpha_n\le {\d {a\over 2}}$
for
$n=1,\,2,\,\ldots\,.$
\item[(ii)]
$\alpha_n$ is nonincreasing and $\alpha_n\ge {\d {a\over 2}}$ for
$n=1,\,2,\,\ldots\,.$
\end{enumerate}
Then
$$\Delta_n(x) \ge 0,\qquad {\rm for}\ \, -a\le x\le a, \quad
n=0,\,1,\,\ldots ,$$
and the equality holds if and only if $n\ge 1$ and $x=\pm a.$

Moreover if (i) is satisfied then
$$\Delta_n(x) < 0,\qquad {\rm for}\ \, |x|> a, \quad n=1,\,2,\,\ldots .$$
\end{theorem}

\noindent{\it Proof.}
By changing  variable $x\to ax$ we can restrict ourselves to the
case $a=1.$
We prove part (i) only, because
the proof of (ii) can be obtained from that of (i) by obvious
modifications.

By assumption we have $p_n(1)=1$ and
$p_n(-1)=(-1)^n.$ Hence $\Delta_n(\pm~1)~=~0.$ Assume now that $|x|<1.$
By (\ref{-x}) it suffices to consider $0\le x<1.$
The proof will go by induction. We have
$\gamma_1\Delta_1(x)=\alpha_1(1-x^2)\ge 0.$ Now assume $\Delta_{n-1}(x) >0.$

In view of $\beta_n=0$ and
$\alpha_n-\alpha_{n-1}=\gamma_{n-1}-\gamma_n,$ Proposition 1
implies
\begin{eqnarray}
\gamma_n\Delta_n& = & \gamma_np_n^2+\alpha_np_{n-1}^2 -
xp_{n-1}p_n,\label{4}\\
\gamma_n\Delta_n& =
&(\alpha_n-\alpha_{n-1})(p_{n-1}^2-p_n^2)+\alpha_{n-1}\Delta_{n-1}.\label{5}
\end{eqnarray}
By (\ref{4}) and the positivity of $x$ we may restrict ourselves to the
case $p_{n-1}(x)p_n(x)~>~0.$
We will assume that $p_{n-1}(x)>0$ and $p_n(x)>0$ (the case
$p_{n-1}(x)<0$ and $p_n(x)<0$ can be dealt with similarly).  By (\ref{5})
and by the induction hypothesis it suffices to consider the case
$p_{n-1}(x)<p_{n}(x),$ since by assumption (i) we have $\alpha_{n-1}\le
\alpha_n.$ In that case since $\gamma_n=1-\alpha_n\ge {1\over
2}\ge \alpha_n$ we get
$$\gamma_np_n(x)-\alpha_np_{n-1}(x)\ge
\alpha_n[p_n(x)-p_{n-1}(x)]\ge 0.$$
Now we apply (\ref{3}) and  obtain
$$\gamma_n\Delta_n \ge (1-x)p_{n-1}(x)p_n(x) >0.$$
The proof of part (i) is thus complete.

We turn to the last part of the statement.
Let (i) be satisfied and $|x|>1.$ By symmetry we can
assume $x>1.$ As before we proceed by induction.
We have
$$\gamma_1\Delta_1(x)= \alpha_1(1-x^2) <1.$$
Assume now that $\Delta_{m}(x)<0$ for $1\le m\le n-1.$ Since
$p_n(1)=1$ and the leading coefficients of $p_n$'s are positive,
the polynomials $p_n$ are positive  for $x>1. $
Thus
$$0>{\Delta_{m}(x)\over p_{m-1}(x)p_m(x)} =
{p_m(x)\over p_{m-1}(x)}- {p_{m+1}(x)\over p_{m}(x)}.$$
for $1\le m\le n-1.$ Hence
$${p_n(x)\over p_{n-1}(x)} \ge \ldots \ge {p_1(x)\over
p_0(x)}=x>1.$$
Now by (\ref{5}) we get
$$\gamma_n\Delta_n \le \alpha_{n-1}\Delta_{n-1}<0.$$
\qed

\noindent
{\bf Remark.} The second part of Theorem 1 is not true under
assumption (ii). Indeed, by (\ref{recsym}), the leading coefficient
of the Tur\'an determinant $\gamma_n\Delta_n(x)$ is equal to
$\gamma_1^{-2}\ldots\gamma_{n-1}^{-2}(\alpha_{n-1}-\alpha_n).$
Thus $\Delta_n(x)$ is positive at infinity  for $n\ge 2.$ One might expect
that in this case $\Delta_n(x)$ is nonnegative on the whole real
axis, but this is not true either. Indeed, it can
be computed that
$$\gamma_1^2\gamma_2\Delta_2(x)=(x^2-1)
[(\gamma_2-\gamma_1)x^2-\alpha_1^2\gamma_2].$$
One can verify that under assumption (ii) we have
$$r:={\alpha_1^2\gamma_2\over \gamma_2-\gamma_1} > 1.$$
(Actually $r\ge 1$ follows from Theorem 1 (ii).)
Hence $\Delta_2(x)<0$ for $1<x<r.$

Sometimes we have to deal with polynomials which are orthogonal
in the interval $[-1,1]$ and normalized at $x=1,$ but the three
term recurrence relation  is not available in  explicit
form. In such cases  the following will be
useful.
\begin{col}
Let the polynomials $p_n$ satisfy
$$xp_n=\gamma_n p_{n+1}+\alpha_np_{n-1},\qquad n=0,\,1,\,\ldots
,$$
with $p_{-1}=0,$  $p_0=1$ and $\alpha_0=0.$
Assume that the sequences $\alpha_n$ and $\alpha_n+\gamma_n$ are
nondecreasing and
$$\lim_{n\to \infty}\alpha_n={1\over 2}a\qquad \lim_{n\to
\infty}\gamma_n={1\over 2}a^{-1},$$
where $0<a<1.$
Then the orthogonality measure for $p_n$ is supported in the
interval $[-1,1].$

Assume that in addition at least one of the following holds
\begin{enumerate}
\item[(i)] $\gamma_n$ is nondecreasing,
\item[(ii)] $\gamma_0\ge 1.$
\end{enumerate}
Then
$$\Delta_n(x)=\widetilde{p}_n^2(x)-\widetilde{p}_{n-1}(x)\widetilde{p}_{n+1}(x)\ge
0  \iff -1\le x\le 1,$$
where
$\widetilde{p}_{n}(x)=p_n(x)/p_n(1).$
\end{col}

\noindent{\it Proof.}
First we will show  that $p_n(1)>0.$ In view of
symmetry of the polynomials this will imply that
the support of the orthogonality measure is contained in
$[-1,1].$

We will show by induction that $p_n(1)/p_{n-1}(1) \ge a >0.$
We have
$${p_0(1)\over p_1(1)}= \gamma_0\le {1\over 2}(a+a^{-1})\le a^{-1}.$$
 Assume that $p_n(1)/p_{n-1}(1) \ge
a.$ Then from recurrence relation we get
$${p_{n+1}(1)\over p_n(1)}= {1\over
\gamma_n}\left (1-\alpha_n{p_{n-1}(1)\over p_n(1)}\right )\ge {1\over
\gamma_n}\left (1-a^{-1}\alpha_n\right ).$$
On the other hand
\begin{eqnarray*}
\gamma_n&=& (\alpha_n+\gamma_n)-\alpha_n\le
{1\over 2}(a+a^{-1})-\alpha_n\\
&\le & {1\over 2}(a+a^{-1})-a^{-2}\alpha_n + (a^{-2}-1)\alpha_n
\\&\le & {1\over 2}(a+a^{-1})-a^{-2}\alpha_n+ (a^{-2}-1){1\over 2}a\\
&=&a^{-1}(1-a^{-1}\alpha_n).
\end{eqnarray*}
Therefore
$${p_{n+1}(1)\over p_n(1)}\ge a.$$

Now we  show  that
$c_n=p_n^2(1)-p_{n-1}(1)p_{n+1}(1)>0$ by induction.
Assume (i).
Similarly to the
proof of Proposition 1 we obtain
\begin{equation}\label{c}
\gamma_nc_n=(\gamma_n-\gamma_{n-1})p_n^2(1)+
(\alpha_n-\alpha_{n-1})p_{n-1}^2 +\alpha_{n-1}c_{n-1}.
\end{equation}
This implies $c_n>0$ for every $n.$

Assume now that (ii) holds and that $c_m>0$ for $m\le n-1.$ Hence
the sequence $p_{m+1}(1)/p_{m}(1)$ is positive and nonincreasing
for $m\le n-1.$ In particular $p_{n}(1)/p_{n-1}(1)\le
p_{1}(1)/p_{0}(1)=\gamma_0^{-1}\le 1.$ Therefore $p_n(1)\le
p_{n-1}(1).$
Rewrite (\ref{c}) in the form
\begin{eqnarray*}
\gamma_nc_n&=& [(\alpha_n+\gamma_n)-(\alpha_{n-1}+\gamma_{n-1})]p_{n}^2(1)\\
&& +
(\alpha_{n}-\alpha_{n-1})[p_{n-1}^2(1)-p_{n}^2(1)]+\alpha_{n-1}c_{n-1}.
\end{eqnarray*}
Thus $c_n>0.$

We have shown that, in both cases (i) and (ii), we have $c_n>0$  and hence
the sequence ${p_{n-1}(1)/ p_n(1)}$ is nondecreasing.
Denote its limit by $r.$ Now plugging  $x=1$ into the recurrence
relation for $p_n,$   dividing both
sides by $p_n(1)$ and taking the limits gives
$$1= {1\over 2}[ar +(ar)^{-1}].$$
Thus $r=a^{-1}.$
Let
$$\widetilde{p}_n(x)={p_n(x)\over p_n(1)}.$$
From the recurrence relation for $p_n$ we obtain
\begin{equation}\label{widetilde}
x\widetilde{p}_n=\widetilde{\gamma}_n\widetilde{p}_{n+1}+
\widetilde{\alpha}_n\widetilde{p}_{n-1} ,
\end{equation}
where
$$\widetilde{\alpha}_n={p_{n-1}(1)\over p_n(1)}\alpha_n,\qquad
\widetilde{\gamma}_n={p_{n+1}(1)\over p_n(1)}\gamma_n.$$
By plugging $x=1$ into (\ref{widetilde}) we get
$$\widetilde{\alpha}_n+\widetilde{\gamma}_n =1.$$
Since both  $p_{n-1}(1)/p_{n}(1)$ and $\alpha_n$ are nondecreasing,
so is $\widetilde{\alpha}_n.$
Moreover it tends to ${1/2}$ at
infinity because the first of its factors tends to $a^{-1}$ while
the second tends to $a/2.$ Therefore the polynomials
$\widetilde{p}_n$ satisfy the assumptions
of Theorem 1 (i). This completes the proof.\qed

\section{The best normalization.}
Assume that the polynomials $p_n$ satisfy (\ref{rec}) and (\ref{norm}).
By multiplying each $p_n$ by a positive constant $\sigma_n$ we obtain
polynomials $p^{(\sigma_n)}_n(x)=\sigma_np_n(x).$ The positivity of
Tur\'an's determinant for the polynomials $p_n$ is not equivalent
to that for the polynomials $p^{(\sigma_n)}_n.$
However, it is possible that the positivity of Turan's determinants in
one normalization implies the positivity in other normalizations.
It turns out that the normalization at  the right most end of the
interval of orthogonality has this feature.
\begin{prop}
Let the polynomials $p_n$ satisfy (\ref{rec}) and (\ref{norm}).
Assume that
$$p_n^2(x)-p_{n-1}(x)p_{n+1}(x)\ge 0,\quad -1\le x\le 1,\quad n\ge 1.$$
Let $p^{(\sigma)}_n(x)=\sigma_np_n(x),$ where $\sigma_n$ is a sequence
of positive constants.  Then
$$\{p^{(\sigma)}_n(x)\}^2-p^{(\sigma)}_{n-1}(x)p^{(\sigma)}_{n+1}(x)\ge
0,\quad-1\le x\le 1,\quad n\ge 1$$
if and only if
$$\sigma_n^2 -\sigma_{n-1}\sigma_{n+1}\ge 0,\quad n\ge 1.$$
\end{prop}

\noindent{\it Proof.}
We have
$$
\{p^{(\sigma)}_n(x)\}^2-p^{(\sigma)}_{n-1}(x)p^{(\sigma)}_{n+1}(x)
$$
$$=
(\sigma_n^2-\sigma_{n-1}\sigma_{n+1})p_n^2(x)+\sigma_{n-1}
\sigma_{n+1}(p_n^2(x)-p_{n-1}(x)p_{n+1}(x)).$$
This shows the "if" part.  On the other hand, since (3) is equivalent to
$p_n(1)=1$ for $n\ge 0,$ we obtain
$$\{p^{(\sigma_n)}_n(1)\}^2-p^{(\sigma_n)}_{n-1}(1)p^{(\sigma_n)}_{n+1}(1)
= \sigma_n^2-\sigma_{n-1}\sigma_{n+1}.$$
This shows the "only if" part.\qed

{\bf Remark.} Proposition 2 says that if the Tur\'an inequality holds
for the polynomials normalized at $x=1$ then it remains true for
any
other normalization if and only if it holds  only at the point $x=1,$
because $p_n^{(\sigma )}(1)=\sigma_n.$

\section{Applications to special symmetric polynomials.}
We will test Theorem 1 on three classes of polynomials:
ultraspherical, q~--~ultraspherical and symmetric Pollaczek
polynomials.

The positivity of Tur\'an's determinants for the first case
is well known (see \cite[p. 209]{erd}).
The ultraspherical polynomials $C_n^{(\lambda)}$ are orthogonal in the interval
$(-1,1)$ with respect to the measure
$(1~-~x^2)^{\lambda~-~(1/2)}~dx,$ where $\lambda >-{1\over 2}.$
When normalized at $x=1$ they satisfy the recurrence relation
$$x\widetilde{C}_n^{(\lambda)}={n+2\lambda\over
2n+2\lambda}\widetilde{C}_{n+1}^{(\lambda)}
+{n\over2n+2\lambda}\widetilde{C}_{n-1}^{(\lambda)} .$$
It can be checked easily that Theorem 1 (i) or (ii) applies
according to $\lambda \ge 0$ or $\lambda \le 0 .$

Let us turn to the $q$--ultraspherical polynomials. They have been
studied by Bustoz and Ismail \cite{bi} but with a normalization
other than the one at the right end of the interval of
orthogonality.  We will exhibit that our  normalization is sharper
in the sense that we can derive the results of
\cite{bi} from ours. Moreover, we will
have no restrictions on the parameters other than that $q$ be
positive.

In standard normalization the $q$--ultraspherical polynomials
are denoted by $C_n(x;\beta|q)$ and they satisfy the recurrence
relation
\begin{equation}\label{ultra}2x C_n(x;\beta|q)= {1-q^{n+1}\over 1-\beta
q^n}C_{n+1}(x;\beta|q)+{1-\beta^2q^{n-1}\over 1-\beta
q^n }C_{n-1}(x;\beta|q).
\end{equation}
The orthogonality
measure is known explicitly (see \cite{ai}, \cite[Thm.  2.2 and Sect.
4]{aw} or \cite[Sect.  7.4]{gr}).  When $|\beta|,\, |q|<1$ it is absolutely
continuous with respect to the Lebesgue measure on the interval $[-1,1].$

\begin{theorem}
Let $0<q<1$ and $|\beta|<1.$ Let $\widetilde{C}_n(x;\beta|q)$ denote
the $q$~--~ultraspherical polynomials normalized at $x=1,$ i.e.
$$\widetilde{C}_n(x;\beta|q)= {{C}_n(x;\beta|q)\over
{C}_n(1;\beta|q)}.$$
Let
$$\Delta_n(x;\beta|q)=\widetilde{C}_n^2(x;\beta|q)-
\widetilde{C}_{n-1}(x;\beta|q)\widetilde{C}_{n+1}(x;\beta|q).$$
Then
$$ \Delta_n(x;\beta|q)\ge 0 \quad \mbox{if and only if}\quad
-1\le x \le 1,$$
with equality only for $x=\pm 1.$
\end{theorem}

\noindent{\it Proof.} The main obstacle in applying Theorem 1 lies in
the fact that the values ${C}_n(1;\beta|q)$ are not given
explicitly. Therefore, we cannot give explicitly the recurrence relation
for $ \widetilde{C}_n(x;\beta|q).$

We will break the proof into two subcases.

\medskip

\noindent(i) $0< \beta <1.$\\
Introduce the polynomials
$$p_n(x)= \beta^{n/2}\prod_{m=1}^n {1-\beta^2 q^{m-1}\over
1-q^m}{C}_n^2(x;\beta|q).$$
Then by (\ref{ultra}) we obtain
$$
xp_n=\gamma_np_{n+1}+ \alpha_np_{n-1},
$$
$$\alpha_n= \beta^{1/2}{1-q^n\over 2(1-\beta q^n)}\qquad
\gamma_n=
\beta^{-1/2}{1-\beta^2q^n\over 2(1-\beta
q^n)}.$$
Observe that
$$ \alpha_n+\gamma_n={1\over 2}(\beta^{1/2}+\beta^{-1/2}).$$
Moreover $\alpha_n $ is nondecreasing and
converges to ${1\over 2}\beta^{1/2}.$ Finally
$$\gamma_0={1\over
2}(\beta^{1/2}+\beta^{-1/2})>1.$$ Therefore we can apply Corollary 1 (ii) with
$a=\beta^{1/2}.$

\medskip

\noindent(ii) $-1 < \beta \le 0.$\\
Introduce the polynomials
$$ p_n(x)= \prod_{m=1}^n {1-\beta^2 q^{m-1}\over
1-q^m}{C}_n^2(x;\beta|q).$$
Then by (\ref{ultra}) we obtain
$$
xp_n=\gamma_np_{n+1}+ \alpha_np_{n-1},
$$
$$\alpha_n= {1-\beta^2q^n\over 2(1-\beta q^n)}\qquad \gamma_n=
{1-q^n\over 2(1-\beta q^n)}.$$
Since both $\alpha_n$ and $\gamma_n$ are increasing sequences
convergent to 1 we can apply Corollary 1(i) with $a=1.$
\qed
\medskip

We turn now to the symmetric Pollaczek polynomials
$P_n^\lambda(x;a).$ They are orthogonal in the interval $[-1,1]$
and satisfy the recurrence relation
$$
xP_n^\lambda(x;a)= {n+1\over 2(n+\lambda+a)}
P_{n+1}^\lambda(x;a) + {n+2\lambda -1\over 2(n+\lambda
+a)}P_{n-1}^\lambda(x;a),$$
where the parameters satisfy $a>0,\ \lambda >0.$
We cannot compute the value $P_n^\lambda(1;a)$ in order to pass
directly to normalization at $x=1.$ Instead, we consider another
auxiliary normalization.
Let
$$p_n(x)= {n!\over (2\lambda)_n}P_n^\lambda(x;a),$$
where $(\mu)_n=\mu(\mu+1)\ldots (\mu +n-1).$
Then the polynomials $p_n$ satisfy the recurrence relation
$$xp_n=
{n+2\lambda\over 2(n+\lambda+a)}p_{n+1}+
{n\over 2(n+\lambda+a)}p_{n-1}.$$
Observe that the  assumptions of
Corollary 1 (i) or (ii) are fulfilled
according to $\lambda\ge a$ or $\lambda\le
a.$ Therefore we have the following.

\begin{theorem}
Let $\lambda>0,\, a>0.$ Let $\widetilde{P}_n^\lambda(x;a)$ denote the
Pollaczek polynomials normalized at $x=1,$ i.e.
$$ \widetilde{P}_n^\lambda(x;a)= {{P}_n^\lambda(x;a)\over
{P}_n^\lambda(1;a)}.$$
Then
$$\{\widetilde{P}_n^\lambda(x;a)\}^2
-\widetilde{P}_{n-1}^\lambda(x;a)\widetilde{P}_{n+1}^\lambda(x;a) \ge 0
\quad \mbox{if and only if}\quad
-1\le x \le 1,$$
with equality only for $x=\pm 1.$
\end{theorem}

\section{Nonsymmetric  polynomials orthogonal in [-1,1]}

In this section we assume that polynomials $p_n$ satisfy
(\ref{rec}) and (\ref{norm})
with $\beta_n$ not necessarily equal to 0 for all $n.$

\begin{theorem}\label{nonsym}
Let  polynomials $p_n$ satisfy  (\ref{rec}) and (\ref{norm}). Let
\begin{equation}\label{init}
|\gamma _0 -\gamma_1|\le \alpha _1\gamma _0-
(\gamma_0-\gamma_1)(1-\gamma_0) .
\end{equation}
Assume that for each $n\ge 2$ one of the following four conditions is
satisfied.
\begin{enumerate}
\item[(i)]
\begin{eqnarray*}
&&\alpha_{n-1}\le \alpha_{n}\le \gamma_n\le \gamma_{n-1},\\
&&{\beta_n+1 + \sqrt{(\beta_n+1)^2-4\alpha_n\gamma_n}\over 2\gamma_n}\le
{ \alpha_{n}-\alpha_{n-1}\over \gamma_{n-1}-\gamma_n}\quad {\rm or} \
(\beta_n+1)^2-4\alpha_n\gamma_n<0.
\end{eqnarray*}
\item[(ii)]
\begin{eqnarray*}
&&\alpha_{n-1}\ge \alpha_{n}\ge \gamma_n\ge \gamma_{n-1},\\
&&{\beta_n+1 - \sqrt{(\beta_n+1)^2-4\alpha_n\gamma_n}\over 2\gamma_n}\ge
{ \alpha_{n-1}-\alpha_{n}\over \gamma_{n}-\gamma_{n-1}}\quad {\rm or} \
(\beta_n+1)^2-4\alpha_n\gamma_n<0.
\end{eqnarray*}
\item[(iii)]
\begin{eqnarray*}
&&\alpha_{n-1}\ge \alpha_{n}\ge {1\over 2},\quad \gamma_{n-1}\ge
\gamma_n\ge {1\over 2},\\
&& { \alpha_{n}-\alpha_{n-1}\over \gamma_{n}-
\gamma_{n-1}}\le {\alpha_n\over \gamma_n}\le 1 \quad {\rm or}
\quad { \alpha_{n}-\alpha_{n-1}\over \gamma_{n}-
\gamma_{n-1}}\ge {\alpha_n\over \gamma_n}\ge 1.
\end{eqnarray*}
\item[(iv)]
\begin{eqnarray*}
&&\alpha_{n-1}\le \alpha_{n},\quad \gamma_{n-1}\le
\gamma_n,\\
&& \cases{\alpha_n\le \gamma_n & \cr
\alpha_n-\alpha_{n-1} \ge \gamma_n-\gamma_{n-1} & \cr} \quad
{\rm or} \quad
\cases{\alpha_n\ge \gamma_n & \cr
\alpha_n-\alpha_{n-1} \le \gamma_n-\gamma_{n-1} & \cr}
\end{eqnarray*} \end{enumerate}
Then
$$\Delta_n(x)=p_n^2(x)-p_{n-1}(x)p_{n+1}(x) \ge 0 \quad {\rm for } \
-1\le x\le 1.$$ \end{theorem}
{\it Proof.} The proof will go by
induction.   Combining (\ref{3}) and (\ref{rec}) for $n=1$ gives
$$\gamma_0^2\gamma_1\Delta_1(x)=(1-x)[(\gamma _0-\gamma_1)(x-\beta
_0)+\alpha _1\gamma_0].  $$ Now using (\ref{norm}) gives that
the positivity of $\Delta_1(x) $ in the interval $[-1,1] $ is equivalent
to (\ref{init}).

Fix $x$ in $[-1,1]$ and assume that
$\Delta_{n}(x)\ge 0. $  Consider two
quadratic functions
\begin{eqnarray*}
A(t)&=&
(t+1)\{(\gamma_n-\gamma_{n-1})t-(\alpha_{n-1}-\alpha_{n})\},\\
B(x;t)&=& \gamma_nt^2 -(\beta_n-x)t + \alpha_n.
\end{eqnarray*}
Set
$$t=-{p_n(x)\over p_{n-1}(x)}.$$
By Proposition 1 it suffices to show that for any $t$
the values $A(t)$ and
$B(x;t)$ cannot be both negative. In order to achieve this we have
to look at the roots of these functions. The roots of $A(t)$ are $-1$ and
$(\alpha_{n-1}-\alpha_{n})/(\gamma_n-\gamma_{n-1});$ hence they
are independent of $x.$ The roots of $B(t)$ have always the same sign
and are equal to
\begin{eqnarray}
r_n^{(1)}(x) & = & {\beta_n
-x-\sqrt{(\beta_n-x)^2-4\alpha_n\gamma_n}\over 2\gamma_n},\label{r1}\\
r_n^{(2)}(x) & = & {\beta_n
-x+\sqrt{(\beta_n-x)^2-4\alpha_n\gamma_n}\over 2\gamma_n}.\label{r2}
\end{eqnarray}
Since the function $u\mapsto u+\sqrt{u^2-a^2},$ $a>0,$ is decreasing for
$u\le -a$ and increasing for $u\ge a$ we have
\begin{eqnarray}\label{roots}
r_n^{(1)}(1) \le r_n^{(1)}(x)\le r_n^{(2)}(x) \le r_n^{(2)}(1)
& \mbox{if} & \beta_n-x\le 0,\\
r_n^{(1)}(-1) \le r_n^{(1)}(x)\le r_n^{(2)}(x) \le r_n^{(2)}(-1)
& \mbox{if} & \beta_n-x\ge 0,
 \end{eqnarray}
provided that $(\beta_n-x)^2-4\alpha_n\gamma_n\ge 0.$
Thus $B(x;t)<0$ implies  $B(1;t)<0$ ($B(-1;t)<0$ respectively) if
$\beta_n-x\le 0$ ($\beta_n-x\ge 0$ respectively).
Hence it suffices to show that the
values $A(t)$ and $B(1;t)$ (the values $A(t)$ and $B(-1;t)$ respectively)
cannot be both negative if $\beta_n-x\le 0$ ($ \beta_n-x\ge 0$
respectively).
We will break the proof into  two  subcases.\medskip\\
{(a)}  $\beta_n -x \le -2\sqrt{\alpha_n\gamma_n}.$

In view of (\ref{3}) and (\ref{1} )the roots of $B(1;t)$ are $-1$ and
$-{\alpha_n\over \gamma_n}.$ By analysing the  positions
of these  numbers with respect to
 the roots of $A(t)$ one can easily verify that under each
of the four assumptions (i) through (iv) the values $A(t)$ and $B(1;t)$
cannot be both
negative.\medskip\\
{(b)}  $\beta_n -x \le -2\sqrt{\alpha_n\gamma_n}.$\\
We  examine the signs of $A(t)$ and $B(-1;t).$
Consider (i), (ii) and (iii). By analysing the mutual position of the
roots of $B(-1;t)$ and $A(t)$ one can verify that $A(t)$ and $B(-1;t)$
cannot be both negative.

In  case (iv) we have that $B(-1;t)\ge 0$ because
$$(1+\beta_n)^2-4\alpha_n\gamma_n=
(2-\alpha_n-\gamma_n)^2-4\alpha_n\gamma_n\le 0.$$
\nopagebreak[4]\qed

{\bf Remark 1.} The assumption (iv) in Theorem 4 does not imply
that the support of the orthogonality measure corresponding to the
polynomials $p_n$ is contained in $[-1,1].$ By (\ref{norm}) we have
$p_n(1)=1$ which implies that the support is located to the left of
$1.$ However, it can extend to the left side beyond $-1.$

{\bf Remark 2.} If we assume that $\beta_n=0$ for $n\ge 0,$ then Theorem
4 reduces to Theorem 1. Indeed, in this case we have
\begin{eqnarray*}
{\beta_n
+1-\sqrt{(\beta_n+1)^2-4\alpha_n\gamma_n}\over 2\gamma_n}&=&\min\left
(1,{\alpha_n\over \gamma_n}
\right ),                  \\
{\beta_n
+1+\sqrt{(\beta_n+1)^2-4\alpha_n\gamma_n}\over 2\gamma_n}&=&\max\left
(1,{\alpha_n\over \gamma_n}
\right ) .
\end{eqnarray*}
\medskip
{\bf Example.}
Set
$$\alpha _{n}= {1\over 2}-{1\over n+2},\quad \gamma _{n}={1\over
2}+{1\over 2(n+2)},\quad \beta_n = {1\over 2(n+2)}. $$
We can check easily that condition (\ref{init}) is satisfied.
We will check that also the assumptions (iii) are satisfied
for every $n\ge 2.$
Clearly we have
 $\alpha_{n-1} \le \alpha_n \le \gamma_{n}\le
\gamma_{n-1}.$ Moreover
\begin{eqnarray*}
r_n^{(2)}(-1)&=&
{\beta_n +1 +\sqrt{(1+\beta_n)^2-4\alpha_n\gamma_n}\over
2\gamma_n}\\
&\le & {\beta_n +1\over \gamma_n}\le 2 = {\alpha_n-\alpha_{n-1}\over
\gamma_{n-1}-\gamma_n}.
\end{eqnarray*}
Let  $p_n(x)$ satisfy (\ref{rec}). By Theorem 4(iii)
$$p_n^2(x)-p_{n-1}(x)p_{n+1}(x)\ge 0 \quad\mbox{for}\quad -1\le x \le 1.
$$
   Let us determine the interval of orthogonality.
    Since $\alpha_n+\beta_n+\gamma_n=1$
we have $p_n(1)=1.  $ Thus the support of the corresponding
orthogonality measure is located to the left of $1.  $ Actually the
support is contained in the interval $[-1,1].$ Indeed, it suffices
to show that $c_n=(-1)^np_n(-1)>0.$ We will show that $c_n\ge
c_{n-1}>0$ by induction. We have $c_0=1.$ Assume $c_n\ge c_{n-1}>0.$
Then by (\ref{rec})
\begin{eqnarray*}
\gamma_nc_{n+1}&=&(1+\beta_n)c_n-\alpha_nc_{n-1}\ge (1+\beta_n-\alpha_n)c_n
 \\
&\ge &
(1-\beta_n-\alpha_n)c_n=\gamma_nc_n.
\end{eqnarray*}
Thus $c_{n+1}\ge c_n >0.$
\section{Polynomials orthogonal in the interval $[0,+\infty).$}
Let $p_n$ be polynomials orthogonal in the positive half axis
normalized at $x=0,$ i.e. $p_n(0)=1.$
Then they satisfy the recurrence relation of the
form
\begin{equation}\label{half}
xp_n= -\gamma_np_{n+1}+(\alpha_n+\gamma_n)p_n-\alpha_{n}p_{n-1},
\quad n=0,\,1,\,\ldots,
\end{equation}
with initial conditions
$p_{-1}=0,$ $p_0=1,$ where $\alpha_n,$
 and $\gamma_n$ are given sequences of real
coefficients such that
\begin{equation}\label{normh}
\alpha_0=0, \  \gamma_0=1, \
\alpha_{n+1}>0,\ \gamma_n >0, \ \mbox{for}\ n=0,\,1,\,\ldots .
\end{equation}
\begin{theorem}
Let polynomials $p_n$ satisfy (\ref{half}) and (\ref{normh}),
and let
$$\alpha_{n-1}\le \alpha_{n},\quad \gamma_{n-1}\le
\gamma_n \quad\mbox{for}\ n\ge 1.$$
Assume that  one of the following two conditions is
satisfied.
\begin{enumerate}
\item[(i)]
$\alpha_n\le \gamma_n \qquad \alpha_n-\alpha_{n-1} \ge
\gamma_n-\gamma_{n-1}.$
\item[(ii)]
$\alpha_n\ge \gamma_n \qquad \alpha_n-\alpha_{n-1} \le \gamma_n-\gamma_{n-1}.$
\end{enumerate}
Then
$$\Delta_n(x)=p_n^2(x)-p_{n-1}(x)p_{n+1}(x) \ge 0 \quad {\rm for } \
x\ge 0.$$
\end{theorem}

{\it Proof.} Let $q_n(x)= p_n(1-x).$ Then by (\ref{half}) we obtain
$$xq_n=
\gamma_nq_{n+1}+(1-\alpha_n-\gamma_n)p_n+\alpha_{n}q_{n-1}.$$
We have  $q_n(1)=1.$
Thus the assumptions (iv) of Theorem \ref{nonsym} are satisfied
for every $n.$ From the proof of Theorem \ref{nonsym} (iv) it follows
that $q_n^2(x)~-~q_{n-1}(x)q_{n+1}(x)~\ge~0$ for $x\le 1$ (the
assumption $x\ge -1$ is inessential).
Taking into account the relation between $p_n$ and $q_n$ gives
the conclusion.
\qed

\medskip

A special case of Theorem 5 is when $\alpha_n-\alpha_{n-1} =
\gamma_n-\gamma_{n-1}$ for every $n.$ In this case, applying (\ref{3})
gives the following.
\begin{prop}
Let polynomials $p_n$ satisfy (\ref{half}) and (\ref{normh}), and let
$$\alpha_n-\alpha_{n-1} =
\gamma_n-\gamma_{n-1}, \qquad n\ge 1.$$
Then
$$p_n^2(x)-p_{n-1}(x)p_{n+1}(x)=\sum_{k=1}^{n}
{(\alpha_k-\alpha_{k-1})\alpha_k\alpha_{k+1}\ldots \alpha_{n-1}\over
\gamma_k\gamma_{k+1}\ldots\gamma_n}(p_k(x)-p_{k-1}(x))^2.$$
In particular, if $\alpha_{n}\ge \alpha_{n-1}$ for  $n\ge 1,$
then
$$p_n^2(x)-p_{n-1}(x)p_{n+1}(x)\ge 0 \quad \mbox{for}\ -\infty
<x<\infty,$$ where equality
holds only for $x=0.$
\end{prop}

\medskip

{\bf Example.}

Let $p_n(x)= L_n^{\alpha}(x)/ L_n^{\alpha}(1),$ where $L_n^{\alpha}(x)$
denote the Laguerre polynomials of order $\alpha >-1.$ Then the
polynomials $p_n$
satisfy
$$xp_n= -(n+\alpha+1)p_{n+1} + (2n+\alpha +1)p_n -np_n.$$
Then
$$\alpha_n-\alpha_{n-1}=\gamma_n-\gamma_{n-1} =1, \quad n\ge 1.$$
Thus Proposition 3 applies.
The formula for $p_n^2-p_{n-1}p_{n+1}$ in this case is not new.
It has been discovered by V. R. Thiruvenkatachar and T. S. Nanjundiah
\cite{tn} (see also \cite[4.7]{asc}.

{\bf Acknowledgement.} I am grateful to J. Bustoz and M. E. H. Ismail for
kindly sending me a preprint of \cite{bi1}.  I thank George Gasper for
pointing out the references  \cite{asc,tn}.

\obeylines{
Institute of Mathematics
Polish Academy of Sciences
ul. Kopernika
00--950 Wroc{\l}aw, Poland
\medskip

{\it Current address:}
Institute of Mathematics
Wroc{\l}aw University
pl. Grunwaldzki 2/4
50--384 Wroc{\l}aw, Poland}

\end{document}